# ANCHORED EXPANSION, PERCOLATION AND SPEED

By Dayue Chen[1] and Yuval Peres[2] with an Appendix by Gábor Pete[3]

*Peking University and University of California, Berkeley*

Benjamini, Lyons and Schramm [*Random Walks and Discrete Potential Theory* (1999) 56–84] considered properties of an infinite graph $G$, and the simple random walk on it, that are preserved by random perturbations. In this paper we solve several problems raised by those authors. The anchored expansion constant is a variant of the Cheeger constant; its positivity implies positive lower speed for the simple random walk, as shown by Virág [*Geom. Funct. Anal.* **10** (2000) 1588–1605]. We prove that if $G$ has a positive anchored expansion constant, then so does every infinite cluster of independent percolation with parameter $p$ sufficiently close to 1; a better estimate for the parameters $p$ where this holds is in the Appendix. We also show that positivity of the anchored expansion constant is preserved under a random stretch if and only if the stretching law has an exponential tail. We then study a simple random walk in the infinite percolation cluster in Cayley graphs of certain amenable groups known as "lamplighter groups." We prove that zero speed for a random walk on a lamplighter group implies zero speed for random walk on an infinite cluster, for any supercritical percolation parameter $p$. For $p$ large enough, we also establish the converse.

**1. Introduction.** Grimmett, Kesten and Zhang (1993) showed that a simple random walk on the infinite cluster of supercritical Bernoulli percolation in $\mathbb{Z}^d$ is transient for $d \geq 3$; in other words, in Euclidean lattices, transience is preserved when the whole lattice is replaced by an infinite percolation cluster. Benjamini, Lyons and Schramm (1999), abbreviated as BLS (1999) hereafter, initiated a systematic study of the properties of a transitive graph

Received March 2003; revised June 2003.
[1]Supported by Ministry of Science and Technology of China Grant G1999075106.
[2]Supported by NSF Grant DMS-01-04073 and by a Miller Professorship at UC Berkeley.
[3]Supported by Hungarian National Foundation for Scientific Research Grant T30074.

$\mathbb{G}$ that are preserved under random perturbations such as passing from $\mathbb{G}$ to an infinite percolation cluster. They conjectured that positivity of the speed for a simple random walk is preserved, and proved this for nonamenable Cayley graphs. Our results (see Theorems 1.5 and 1.6) lend further support to this conjecture.

We first consider the stability of a related geometric quantity. Denote by $V(\mathbb{G})$ and $E(\mathbb{G})$, respectively, the sets of vertices and edges of an infinite graph $\mathbb{G}$. For $S \subset V(\mathbb{G})$, denote by $|S|$ the cardinality of $S$ and by $\partial S = \partial_{\mathbb{G}} S$ the set of edges that have one end in $S$ and the other in $S^c$. We say that $S$ is *connected* if the induced subgraph on $S$ is connected. Fix $o \in V(\mathbb{G})$. The *anchored expansion constant* of $\mathbb{G}$,

$$\imath_E^*(\mathbb{G}) := \lim_{n \to \infty} \inf \left\{ \frac{|\partial S|}{|S|} : o \in S \subset V(\mathbb{G}), S \text{ is connected}, n \leq |S| < \infty \right\}$$

was defined in BLS (1999). The quantity $\imath_E^*(\mathbb{G})$ does not depend on the choice of the basepoint $o$. It is related to the isoperimetric constant

$$\imath_E(\mathbb{G}) := \inf \left\{ \frac{|\partial S|}{|S|} : S \subset V(\mathbb{G}), S \text{ is connected}, 1 \leq |S| < \infty \right\},$$

but as we shall see, $\imath_E^*(\cdot)$ is more robust. BLS (1999) asked if the positivity of $\imath_E^*(\mathbb{G})$ is preserved when $\mathbb{G}$ undergoes a random perturbation.

In *p-Bernoulli bond percolation* in $\mathbb{G}$, each edge of $\mathbb{G}$ is independently declared *open* with probability $p$ and *closed* with probability $1 - p$. Thus a bond percolation $\omega$ is a random subset of $E(\mathbb{G})$. We usually identify the percolation $\omega$ with the subgraph of $\mathbb{G}$ consisting of all open edges and their end-vertices. A connected component of this subgraph is called an *open cluster*, or simply a *cluster*. The probability that there is an infinite cluster is monotone in $p$. Let $p_c = p_c(\mathbb{G}) = \inf\{p: \text{there is an infinite cluster a.s.}\}$. When $p \in (p_c, 1)$, with positive probability the open cluster $\mathbb{H}$ that contains $o$ is infinite; it is easy to see that $\imath_E(\mathbb{H}) = 0$ a.s.

Theorem 2 of Benjamini and Schramm (1996) states that $p_c(\mathbb{G}) \leq 1/(\imath_E(\mathbb{G}) + 1)$, but their proof yields the stronger inequality $p_c(\mathbb{G}) \leq 1/(\imath_E^*(\mathbb{G}) + 1)$.

THEOREM 1.1. *Consider p-Bernoulli percolation on a graph $G$ with $\imath_E^*(G) > 0$. If $p < 1$ is sufficiently close to $1$, then almost surely on the event that the open cluster $\mathbb{H}$ containing $o$ is infinite, we have $\imath_E^*(\mathbb{H}) > 0$.*

Our proof of Theorem 1.1, given in the next section, shows the conclusion holds for all $p > 1 - h/(1+h)^{1+1/h}$, where $h = \imath_E^*(\mathbb{G})$. A refinement of the argument, due to Gábor Pete (see the Appendix) shows the conclusion holds for all $p > 1/(\imath_E^*(\mathbb{G}) + 1)$. The Appendix also contains the analog of Theorem 1.1 for site percolation.

Next, let $\mathbb{G}$ be an infinite graph of bounded degree and pick a probability distribution $\nu$ on the positive integers. Replace each edge $e \in E(\mathbb{G})$ by a path that consists of $L_e$ new edges, where the random variables $\{L_e\}_{e \in E(\mathbb{G})}$ are independent with law $\nu$. Let $\mathbb{G}^\nu$ denote the random graph obtained in this way. We call $\mathbb{G}^\nu$ a *random stretch* of $\mathbb{G}$. If the support of $\nu$ is unbounded, then $\imath_E(\mathbb{G}^\nu) = 0$ a.s. Say that $\nu$ has an *exponential tail* if $\nu[\ell, \infty) < e^{-\varepsilon \ell}$ for some $\varepsilon > 0$ and all sufficiently large $\ell$.

THEOREM 1.2. *Suppose that $\mathbb{G}$ is an infinite graph of bounded degree and $\imath_E^*(\mathbb{G}) > 0$. If $\nu$ has an exponential tail, then $\imath_E^*(\mathbb{G}^\nu) > 0$ a.s.*

On the other hand, if $\nu$ has a tail that decays slower than exponentially, then taking the binary tree as $\mathbb{G}$, we have $\imath_E^*(\mathbb{G}) > 0$ yet $\imath_E^*(\mathbb{G}^\nu) = 0$ a.s. See Remark 2.2 in the next section.

By a *Galton–Watson tree* we mean a family tree of a Galton–Watson process.

COROLLARY 1.3. *For a supercritical Galton–Watson tree $\mathbb{T}$, given nonextinction, we have $\imath_E^*(\mathbb{T}) > 0$ a.s.*

Theorem 1.2 and Corollary 1.3 answer Questions 6.3 and 6.4 of BLS (1999), while Theorem 1.1 partially answers Question 6.5 of the same paper.

The importance of anchored expansion is exhibited by the following theorem, conjectured in BLS [(1999), Conjecture 6.2]. For a vertex $x$, denote by $|x| = |x|_\mathbb{G}$ the distance (the least number of edges on a path) from $x$ to the basepoint $o$ in $\mathbb{G}$.

THEOREM 1.4 [Virág (2000)]. *Let $\mathbb{G}$ be a bounded degree graph with $\imath_E^*(\mathbb{G}) > 0$. Then the simple random walk $\{X_n\}$ in $\mathbb{G}$, started at $o$, satisfies $\liminf_{n \to \infty} |X_n|/n > 0$ a.s. and there exists $C > 0$ such that $\mathbb{P}[X_n = o] \leq \exp(-Cn^{1/3})$ for all $n \geq 1$.*

Earlier, Thomassen (1992) showed that a condition weaker than $\imath_E^*(\mathbb{G}) > 0$ suffices for transience of the random walk $\{X_n\}$. As noted in Virág (2000), Theorem 1.4, in conjunction with Corollary 1.3, implies that the speed of simple random walk on supercritical Galton–Watson trees is positive, a result first proved in Lyons, Pemantle and Peres (1995). Other applications of anchored expansion are in Häggström, Schonmann and Steif (2000).

In Section 3 we address another problem in BLS (1999) concerning the *speed* $\lim_n \frac{|X_n|}{n}$ of a random walk $\{X_n\}$. Consider again $p$-Bernoulli bond percolation in $\mathbb{G}$ with parameter $p > p_c(\mathbb{G})$. Theorem 1.3 of BLS (1999) states that if the Cayley graph $\mathbb{G}$ is nonamenable [i.e., $\imath_E(G) > 0$], then a simple random walk in an infinite cluster of Bernoulli percolation on $\mathbb{G}$ has positive speed. On the other hand, if a graph $\mathbb{G}$ has subexponential growth, that is, if $\limsup |\{x \in V(\mathbb{G}) : |x| \leq n\}|^{1/n} = 1$, then a simple random walk on $\mathbb{G}$ (and on any subgraph) has zero speed [Varopoulos (1985)]. It is therefore natural to study, as suggested in BLS (1999), a simple random walk in the infinite cluster of an amenable Cayley graph with exponential growth.

The lamplighter groups $G_d$ are amenable groups with exponential growth, introduced by Kaimanovich and Vershik (1983). The corresponding Cayley graphs $\mathbb{G}_d$ (for the standard generators) can be described as follows. A vertex of $\mathbb{G}_d$ can be identified as $(m, \eta) \in \mathbb{Z}^d \times \{\text{finite subsets of } \mathbb{Z}^d\}$. Heuristically, $\mathbb{Z}^d$ is the set of lamps, $\eta$ is the set of lamps which are on, and $m$ is the position of the lamplighter, or "marker." In each step, either the lamplighter switches the current lamp (from on to off, or from off to on) or moves to one of the neighboring sites in $\mathbb{Z}^d$. Each vertex in $\mathbb{G}_d$ has degree $2d + 1$; one edge corresponds to flipping the state of the lamp at location $m$, and the other $2d$ edges correspond to moving the marker. For example, if $d = 1$, the neighbors of $(m, \eta)$ are $(m + 1, \eta)$, $(m - 1, \eta)$ and $(m, \eta \Delta \{m\})$, where $\eta \Delta \{m\}$ is $\eta \setminus \{m\}$ if $m \in \eta$, and is $\eta \cup \{m\}$ if $m \notin \eta$. For a more detailed description see Lyons, Pemantle and Peres (1996). Kaimanovich and Vershik (1983) showed that simple random walk in $\mathbb{G}_d$ has speed zero for $d = 1, 2$ and has positive speed for $d \geq 3$.

We now study simple random walk $\{X_n\}$ in the unique infinite cluster of $p$-Bernoulli bond percolation in $\mathbb{G}_d$. If $x$ is a vertex in the open cluster containing $o$, let $|x|_\omega$ be the graph distance in this cluster from $x$ to $o$.

THEOREM 1.5.   *Let $d \in \{1, 2\}$. Then the simple random walk in the infinite cluster of $\mathbb{G}_d$ has zero speed, that is, $\lim_n \frac{|X_n|_\omega}{n} = 0$ a.s. on the event that $o$ is in the infinite cluster.*

THEOREM 1.6.   *Suppose that $d \geq 3$. If $p > p_c(\mathbb{Z}^d)$, then the simple random walk in the infinite cluster of Bernoulli bond percolation in $\mathbb{G}_d$ has positive speed.*

These results support Conjectures 1.4 and 1.5 of BLS (1999); they are extended in Theorems 3.1 and 3.2.

**2. Anchored expansion.**   The idea of the following lemma is from Kesten (1982).

LEMMA 2.1.   *Let $\mathcal{A}_n = \{S \subset V(\mathbb{G}) : o \in S, \ S \text{ is connected}, \ |\partial S| = n\}$. If $i_E^*(\mathbb{G}) > h > 0$, then for all sufficiently large $n$,*

$$|\mathcal{A}_n| \leq [\Psi(h)]^n,$$

*where*

$$\Psi(h) = (1+h)^{1+1/h}/h.$$

PROOF. Consider $p$-Bernoulli bond percolation in $\mathbb{G}$. Let $\mathbb{H}$ be the open cluster containing $o$. Then $V(\mathbb{H})$ is the set of vertices which can be reached from $o$ via open bonds. For any $S \in \mathcal{A}_n$, a spanning tree on $S$ has $|S| - 1$ edges. Also, note that $|\partial S| \geq h|S|$ if $n = |\partial S|$ is large enough. Therefore

$$P(V(\mathbb{H}) = S) \geq p^{|S|-1}(1-p)^{|\partial S|} \geq p^{n/h-1}(1-p)^n,$$

whence

$$1 \geq P(V(\mathbb{H}) \in \mathcal{A}_n) = \sum_{S \in \mathcal{A}_n} P(V(\mathbb{H}) = S) \geq |\mathcal{A}_n| p^{n/h-1}(1-p)^n.$$

Thus, if $n$ is sufficiently large,

$$|\mathcal{A}_n| \leq \left(\frac{1}{p}\right)^{n/h-1}\left(\frac{1}{1-p}\right)^n$$

holds for any $p \in (0,1)$. Letting $p = 1/(1+h)$ concludes the proof. $\square$

PROOF OF THEOREM 1.1. Let $\mathbb{H}$ be the open cluster containing $o$. Denote

$$\mathcal{A}_n(\mathbb{H}) = \{S \subset V(\mathbb{H}) : o \in S, S \text{ is connected in } \mathbb{H} \text{ and } |\partial_\mathbb{G} S| = n\}.$$

Suppose that $S \in \mathcal{A}_n(\mathbb{H})$. Then $S$ is also a connected subset of $V(\mathbb{G})$; we shall use a subscript to indicate the graph considered. For $S \in \mathcal{A}_n$, each edge in $\partial_\mathbb{G} S$ is independently open with probability $p$. By the large deviation principle [see Dembo and Zeitouni (1998), Theorem 2.1.14],

$$(2.1) \qquad P\left(S \in \mathcal{A}_n(\mathbb{H}), \frac{|\partial_\mathbb{H} S|}{|\partial_\mathbb{G} S|} \leq \alpha\right) \leq e^{-nI_p(\alpha)},$$

where the rate function $I_p(\alpha) = \alpha \log \frac{\alpha}{p} + (1-\alpha)\log \frac{1-\alpha}{1-p}$ satisfies $I_p(\alpha) > 0$ for $\alpha < p$. Recall $\Psi(h)$ defined in Lemma 2.1. When $p > 1 - 1/\Psi(h)$, we have $I_p(0) = -\log(1-p) > \log \Psi(h)$, so there exists $\alpha_0 > 0$ such that $I_p(\alpha_0) > \log \Psi(h)$. For sufficiently large $n$,

$$P\left(\exists S \in \mathcal{A}_n(\mathbb{H}) : \frac{|\partial_\mathbb{H} S|}{|\partial_\mathbb{G} S|} \leq \alpha_0\right) \leq |\mathcal{A}_n| e^{-nI_p(\alpha_0)} \leq \Psi(h)^n e^{-nI_p(\alpha_0)},$$

which is summable in $n$. By the Borel–Cantelli lemma,

$$\liminf_{n\to\infty} \left\{\frac{|\partial_\mathbb{H} S|}{|\partial_\mathbb{G} S|} : o \in S \subset V(\mathbb{H}), S \text{ is connected}, n \leq |\partial_\mathbb{G} S|\right\} \geq \alpha_0 \qquad \text{a.s.,}$$

6                     D. CHEN AND Y. PERESwhence

$$\liminf_{n \to \infty} \left\{ \frac{|\partial_{\mathbb{H}} S|}{|S|} : o \in S \subset V(\mathbb{H}), \right.$$

$$\left. S \text{ is connected}, n \leq |S| < \infty \right\} \geq \alpha_0 \iota_E^*(\mathbb{G}) \qquad \text{a.s.} \qquad \square$$

PROOF OF THEOREM 1.2. Let $L_1, L_2, \ldots, L_n$ be i.i.d. random variables with distribution $\nu$. Since $\nu$ has an exponential tail, there is an increasing convex rate function $I(\cdot)$ such that $I(c) > 0$ for $c > EL_i$ and $P(\sum_{i=1}^{n} L_i > cn) \leq \exp(-nI(c))$ for all $n$ [see Dembo and Zeitouni (1998), Theorem 2.2.3, page 27]. Choose $c$ large enough such that $I(c) > \log \Psi(h)$. For any $S \in \mathcal{A}_n$, let $\mathsf{Edge}(S)$ be the set of edges with at least one end in $S$. Note that $|\partial S| \leq |\mathsf{Edge}(S)| \leq D|S|$, where $D$ is the maximal degree in $\mathbb{G}$. Thus for $S \in \mathcal{A}_n$,

$$P\left(\frac{\sum_{e \in \mathsf{Edge}(S)} L_e}{D|S|} > c\right) \leq P\left(\frac{\sum_{i=1}^{D|S|} L_i}{D|S|} > c\right)$$
$$\leq \exp(-D|S|I(c)) \leq \exp(-|\partial S|I(c)).$$

Therefore for all $n$,

$$P\left(\exists S \in \mathcal{A}_n : \frac{\sum_{e \in \mathsf{Edge}(S)} L_e}{D|S|} > c\right) \leq |\mathcal{A}_n| e^{-I(c)n},$$

which is summable. By the Borel–Cantelli lemma, with probability 1, for any sequence of sets $\{S_n\}$ such that $S_n \in \mathcal{A}_n$ for each $n$, we have

$$\limsup_{n \to \infty} \frac{\sum_{e \in \mathsf{Edge}(S_n)} L_e}{D|S_n|} \leq c \qquad \text{a.s.}$$

Therefore

$$\liminf_{n \to \infty} \left\{ \frac{|\partial S|}{\sum_{e \in \mathsf{Edge}(S)} L_e} : o \in S \subset V(\mathbb{G}), \right.$$

$$\left. S \text{ is connected}, n \leq |\partial S| \right\} \geq \frac{h}{Dc} \qquad \text{a.s.}$$

Since $\mathbb{G}^\nu$ is obtained from $\mathbb{G}$ by adding new vertices, $V(\mathbb{G})$ can be embedded into $V(\mathbb{G}^\nu)$ as a subset. In particular, we can choose the same basepoint $o$ in $\mathbb{G}^\nu$ and in $\mathbb{G}$. For $S$ connected in $\mathbb{G}$ such that $o \in S \subset V(\mathbb{G})$, there is a unique *maximal* connected $\widetilde{S} \subset V(\mathbb{G}^\nu)$ such that $\widetilde{S} \cap V(\mathbb{G}) = S$; it satisfies $|\widetilde{S}| \leq \sum_{e \in \mathsf{Edge}(S)} L_e$. In computing $\iota_E^*(\mathbb{G}^\nu)$ it suffices to consider only such maximal $\widetilde{S}$'s, so we conclude that $\iota_E^*(\mathbb{G}^\nu) \geq h/dc > 0$. $\square$

REMARK 2.2. Suppose that the distribution $\nu$ of $L$ does not have an exponential tail. Then for any $c > 0$ and any $\varepsilon > 0$, we have $P(\sum_{i=1}^n L_i \geq cn) \geq P(L_1 \geq cn) \geq e^{-\varepsilon n}$ for infinitely many $n$'s, where $\{L_i\}$ are i.i.d. with law $\nu$. Let $\mathbb{G}$ be a binary tree with the root $o$ as the basepoint. Pick a collection of $2^n$ pairwise disjoint paths from level $n$ to level $2n$:

$$P\left(\text{along at least one of these } 2^n \text{ paths } \sum_{i=1}^n L_i \geq cn\right)$$
$$\geq 1 - (1 - e^{-\varepsilon n})^{2^n} \geq 1 - \exp(-e^{-\varepsilon n} 2^n) \to 1.$$

With probability very close to 1 (depending on $n$), there is a path from level $n$ to $2n$ along which $\sum_{i=1}^n L_i \geq cn$. Take such a path and extend it to the root $o$. Let $S$ be the set of vertices in the extended path from the root $o$ to level $2n$. Then

$$\frac{|\partial S|}{\sum_{e \in \mathsf{Edge}(S)} L_e} \leq \frac{2n+1}{cn} \approx \frac{2}{c}.$$

Since $c$ can be arbitrarily large, $\iota_E^*(\mathbb{G}^\nu) = 0$ a.s. This shows that the exponential tail condition is necessary to ensure the positivity of $\iota_E^*(\mathbb{G}^\nu)$.

PROOF OF COROLLARY 1.3. A Galton–Watson process is uniquely determined by the offspring distribution $\{p_0, p_1, p_2, \ldots\}$. Let $\mathbb{T}$ be a Galton–Watson tree and $o$ its root.

*Case* (i) $p_0 = p_1 = 0$. For any finite $S \subset V(\mathbb{T})$, $|S| \leq |\partial S|(\frac{1}{2} + \frac{1}{2^2} + \cdots) \leq |\partial S|$. So $\iota_E^*(\mathbb{T}) \geq \iota_E(\mathbb{T}) \geq 1$.

*Case* (ii) $p_0 = 0, p_1 > 0$. In this case the Galton–Watson tree $\mathbb{T}$ can be viewed as random stretch $\mathbb{G}^\nu$ of another Galton–Watson tree $\mathbb{G}$, where $\mathbb{G}$ is generated according to $p'_k = p_k/(1-p_1)$, $k = 2, 3, \ldots$, $p'_0 = p'_1 = 0$ and $\nu$ is the geometric distribution with parameter $p_1$. By Theorem 1.2, $\iota_E^*(\mathbb{T}) = \iota_E^*(\mathbb{G}^\nu) > 0$ a.s.

*Case* (iii) $p_0 > 0$. Let $f(s) = \sum_{i=0}^\infty p_k s^k$ and let $q < 1$ be the extinction probability, so that $q = f(q)$. An infinite Galton–Watson tree can be constructed as follows; see Lyons (1992). Begin with the root which is declared to be *open*. Add to the root a random number of edges according to probability distribution $P(Y = k) = p_k(1 - q^k)/(1 - q)$. Declare each vertex *open* with probability $1 - q$ and *closed* with probability $q$, independently of each other. If all the newly added vertices are closed, discard the entire assignment and reassign open/closed all over again. For each open vertex, repeat the same procedure. For each closed vertex, attach to it independently a Galton–Watson tree conditioned to be finite.

The subtree $\mathbb{T}_1$ consisting of open vertices and edges connecting them is a Galton–Watson tree without leaves, and $i_E^*(\mathbb{T}_1) > 0$ according to case (ii). For each open vertex $x$ of $\mathbb{G}$, label its offspring from 1 to $Y_x$, where $Y_x$ is a random variable with $P(Y_x = k) = p_k(1 - q^k)/(1 - q)$. Along the sequence of $Y_x$ vertices, each is open with probability $1 - q$ and closed with probability $q$ (independently of each other if we ignore the constraint that there is at least one open vertex). The number of closed vertices before the first open vertex is stochastically bounded above by a random variable with a geometric distribution. The same statement holds for the number of closed vertices after the last open vertex, and for the number of closed vertices between the $k$th open vertex and the $(k+1)$st open vertex.

Let $L_1$ be the total number of vertices of finite Galton–Watson trees attached to the closed vertices before the second open vertex (if it ever exists). Similarly, let $L_2$ be the total number of vertices of finite Galton–Watson trees attached to the closed vertices between the second open vertex and the third open vertex (if it ever exists). And so on, until the last open vertex among the offspring of $x$. The variables $L_2, L_3, \ldots$ are i.i.d.; $L_1$ is independent of other $L_i$'s but has a different distribution. Thus we may identify the Galton–Watson tree $\mathbb{T}$ as a random stretch of $\mathbb{T}_1$ in computing $i_E^*(\mathbb{T})$. Although there are two different distributions in the random stretch, the same argument works since both have exponential tails.

All $L_i$'s are stochastically dominated by $\sum_{j=1}^{W_1+W_2} U_j$, where $W_1, W_2, U_1, U_2, \ldots$ are random variables, independent of each other, $P(W_i = k) = q^k(1 - q)$, $k = 0, 1, 2, \ldots$, and $U_j$ is the size of a Galton–Watson tree conditioned on extinction. Let $\nu$ be the probability distribution of $\sum_{j=1}^{W_1+W_2} U_j$. By the next lemma we conclude that $\nu$ has an exponential tail. Applying Theorem 1.2 completes the proof.
□

LEMMA 2.3 [Harris (1963), Theorem 13.1]. *For a supercritical Galton–Watson process, the size of a Galton–Watson tree conditioned on extinction has a distribution with an exponential tail.*

**3. Speed of random walk.** We start with a generalization of the lamplighter groups defined in the Introduction.

Let $\mathbb{G}$ be the Cayley graph of a finitely generated infinite group $\mathsf{G}$ with a fixed, symmetric (i.e., closed under inversion) set of generators. We identify vertices of $\mathbb{G}$ with elements of the group $\mathsf{G}$. Two points $x$ and $y$ of $\mathbb{G}$ are neighbors if $xy^{-1}$ is a generator.

Let $\mathbb{F}$ be the Cayley graph of a finite group $\mathsf{F}$ generated by a fixed symmetric set of generators.

By $\sum_{x \in \mathbb{G}} \mathbb{F}$ we denote the set of elements of $\mathbb{F}^{\mathbb{G}}$ such that at most finitely many of the coordinates are not the identity element of $\mathsf{F}$. An element of $\sum_{x \in \mathbb{G}} \mathbb{F}$ is called a *configuration* and is denoted by $\eta = \{\eta(x) : x \in V(\mathbb{G})\}$, where $\eta(x) \in V(\mathbb{F})$ is the $x$-coordinate of $\eta$. We will sometimes write $x \in \mathbb{G}$ as an abbreviation for $x \in V(\mathbb{G})$.

Define a new graph $\mathbb{W} = \mathbb{G} \ltimes \sum_{x \in \mathbb{G}} \mathbb{F}$ as a semidirect product of $\mathbb{G}$ with the direct sum of copies of $\mathbb{F}$ indexed by $\mathbb{G}$. Vertices of $\mathbb{W}$ are identified as $\{(m, \eta) : m \in V(\mathbb{G}), \eta \in \sum_{x \in \mathbb{G}} \mathbb{F}\}$. Two vertices, $(m, \eta)$ and $(m_1, \xi)$, are neighbors if either:

(i) $m = m_1$, $\eta(x) = \xi(x)$ for all $x \neq m$, and $\eta(m)$ is a neighbor of $\xi(m)$ in $\mathbb{F}$, or

(ii) $\eta = \xi$, and $m, m_1$ are neighbors in $\mathbb{G}$.

In particular, if $\mathsf{F} = \{0, 1\}$ is the group of two elements and $\mathbb{G}$ is $\mathbb{Z}^d$, then $\mathbb{G} \ltimes \sum_{x \in \mathbb{G}} \mathbb{F}$ is exactly $\mathbb{G}_d$ described before Theorem 1.5. Also note that the above definition applies to the case of an arbitrary graph $\mathbb{G}$ and a finite group $\mathsf{F}$, as well.

From now on $\mathbb{G}$ will be an infinite amenable Cayley graph. Then the graph $\mathbb{G} \ltimes \sum_{x \in \mathbb{G}} \mathbb{F}$ is amenable and grows exponentially. By Burton and Keane (1989), there is only one infinite cluster when percolation occurs.

We say that $\mathbb{G}$ is *recurrent* if the simple random walk in $\mathbb{G}$ is recurrent; this is equivalent to $\mathsf{G}$ being a finite extension of $\mathbb{Z}^1$ or $\mathbb{Z}^2$ [see, e.g., Woess (2000), Theorem 3.24, page 36]. The following theorem is a generalization of Theorem 1.5.

THEOREM 3.1. *Suppose that $\mathbb{G}$ is a recurrent Cayley graph and that $\mathbb{F}$ is the Cayley graph of a finite group. Then the simple random walk in the infinite cluster of supercritical Bernoulli bond percolation in $\mathbb{W} = \mathbb{G} \ltimes \sum_{x \in \mathbb{G}} \mathbb{F}$ has zero speed a.s.*

On the other hand, if $\mathbb{G}$ is a transient Cayley graph of polynomial or exponential growth, then for $p$ sufficiently close to 1, the infinite cluster of $p$-Bernoulli bond percolation in $\mathbb{G}$ is transient. For $\mathbb{G} = \mathbb{Z}^d, d \geq 3$ and any $p > p_c(\mathbb{G})$, this is due to Grimmett, Kesten and Zhang (1993); for other Cayley graphs of polynomial growth it is due to Benjamini and Schramm (1998); see also Theorem 9 in Angel, Benjamini, Berger and Peres (2004); for Cayley graphs of exponential growth it is Theorem 1.8 of BLS (1999). The following theorem generalizes our Theorem 1.6.

THEOREM 3.2. *Let $0 < p < 1$. Suppose that the infinite cluster of $p$-Bernoulli bond percolation in the Cayley graph $\mathbb{G}$ is transient and that $\mathbb{F}$ is the Cayley graph of a finite group. Then the simple random walk in the infinite cluster of $p$-Bernoulli bond percolation in $\mathbb{W} = \mathbb{G} \ltimes \sum_{x \in \mathbb{G}} \mathbb{F}$ has positive speed a.s.*

Fix a vertex $o$ of $\mathbb{W} = \mathbb{G} \ltimes \sum_{x \in \mathbb{G}} \mathbb{F}$ as the basepoint, for example, the vertex corresponding to the unit element of the group. Let $\|x\|$ be the distance between the vertex $x$ and the basepoint $o$ in $\mathbb{W}$. Certainly, $\|x\| \leq |x|_\omega$. In the other direction, Lemma 4.6 of BLS (1999) states that if $\lim_n \|X_n\|/n = 0$, then $\lim_n |X_n|_\omega/n = 0$. For this reason we shall consider $\|x\|$ instead of $|x|_\omega$.

It will be useful to consider *delayed simple random walk* $Z = Z^\omega$ on $\omega$, defined as follows. Let $Z(0)$ be some fixed vertex of $\mathbb{W} = G \ltimes \sum_{x \in \mathbb{G}} \mathbb{F}$. For $n \geq 0$, given $\langle Z(0), \ldots, Z(n) \rangle$ and $\omega$, let $Z'(n+1)$ be a uniform random choice from $Z(n)$ and its neighbors in $\mathbb{W}$. Set $Z(n+1) := Z'(n+1)$ if the edge $[Z(n), Z'(n+1)]$ belongs to $\omega$; otherwise, let $Z(n+1) := Z(n)$. By Lemma 4.2 of BLS (1999), the speed $\lim_{n \to \infty} \|Z(n)\|/n$ exists and is constant a.s.

LEMMA 3.3.

$$(3.1) \quad \lim_{n \to \infty} \frac{\|X_n\|}{n} \geq \lim_{n \to \infty} \frac{\|Z(n)\|}{n} \geq c \lim_{n \to \infty} \frac{\|X_n\|}{n} \quad a.s.,$$

*where $c > 0$ is a deterministic constant.*

PROOF. A sample path of $Z$ is obtained from a sample path of $X$ by repeating $X_n$ a random number of times, with a geometric distribution. The parameter of the geometric distribution is in $[1/(D+1), D/(D+1)]$, where $D$ is the degree of a vertex of $\mathbb{W}$. Therefore (3.1) holds. □

$Z$ will always denote the delayed random walk in a cluster $\omega$ in $\mathbb{W}$. Denote by $P_\omega$ the law of $Z$ for fixed $\omega$, and let $E_\omega$ be the corresponding expectation operator. Denote by $\mathbf{E}$ the average over realizations of $\omega$. Write $Z(n) = (m_n, \eta_n)$ and call the first component $m_n$ the *marker*. By Lemma 3.3 and the discussion preceding it, it is enough to determine if $\lim_{n \to \infty} \|Z(n)\|/n$ is positive or not.

Our first goal is to prove Theorem 3.1. A key fact is that the motion of the marker is recurrent in the following sense.

LEMMA 3.4. *Suppose that $\mathbb{G}$ is a recurrent Cayley graph and that $\mathbb{F}$ is the Cayley graph of a finite group. Let $Z$ be as above and write $Z(n) = (m_n, \eta_n)$. Then*

$$(3.2) \quad \mathbf{E} P_\omega(m_n = m_0 \text{ for some } n \geq 1) = 1.$$

PROOF. Introduce the stopping times

$$\tau_N = \min\{n \geq 0; |m_n|_\mathbb{G} = N\},$$
$$\tau_o^+ = \min\{n \geq 1; m_n = o\}.$$

Then (3.2) can be rewritten as

$$\lim_{N\to\infty} \mathbf{E}P_\omega(\tau_N < \tau_o^+ | m_0 = o) = 0.$$

Let $\mathbb{G}_N = \{x \in \mathbb{G} : |x|_\mathbb{G} \leq N\}$. There may be several disjoint clusters in a realization of $p$-Bernoulli bond percolation in the finite graph $\mathbb{G}_N \ltimes \sum_{x \in \mathbb{G}_N} \mathbb{F}$, and each cluster may have several vertices with the marker at $o$. Consider a cluster $\mathbb{H}$ with at least one vertex whose marker is at $o$. If there are $k$ vertices in $\mathbb{H}$ with the marker at $o$, "glue" these $k$ vertices together as one vertex denoted by $\Theta$. Let $\mathbb{H}'$ be the modified graph of the cluster $\mathbb{H}$. Coupling the delayed simple random walks in $V(\mathbb{H})$ and in $V(\mathbb{H}')$, we find that

$$(3.3) \quad \frac{1}{k}\sum_{x \in V(\mathbb{H}), m(x)=o} P_\omega(\tau_N < \tau_o^+ | Z(0) = x) = P_\omega(\tilde{\tau}_N < \tilde{\tau}_o^+ | Z'(0) = \Theta),$$

where

$$\tilde{\tau}_N = \min\{n \geq 0; |m_n|_\mathbb{G} = N\} \quad \text{and} \quad \tilde{\tau}_o^+ = \min\{n \geq 1; m_n = o\}$$

are the stopping times for the delayed simple random walk $Z'$ in $V(\mathbb{H}')$.

The delayed simple random walk $Z'$ in $V(\mathbb{H}')$ is a reversible Markov chain with respect to the measure $\pi$ where $\pi(\Theta) = k$, and $\pi(x) = 1$ for all other $x \in V(\mathbb{H}')$, $x \neq \Theta$. Let $D$ denote the degree of a vertex in $\mathbb{W} = \mathbb{G} \ltimes \sum_{x \in \mathbb{G}} \mathbb{F}$, and let $\mathcal{F}$ be the class of functions with the following properties:

$$(3.4) \quad f : V(\mathbb{H}) \cup \{\Theta\} \to [0,1], \qquad f(x) = \begin{cases} 0, & \text{if } x = \Theta \text{ or } m(x) = o, \\ 1, & \text{if } |m(x)|_\mathbb{G} = N. \end{cases}$$

Applying the Dirichlet principle [see Liggett (1985), page 99], we have that

$$(3.5) \quad \begin{aligned} 2\pi(\Theta) & P_\omega(\tilde{\tau}_N < \tilde{\tau}_o^+ | Z(0) = \Theta) \\ &= \inf_{f \in \mathcal{F}} \sum_{u \in V(\mathbb{H}')} \sum_{[u,v] \in E(\mathbb{H}')} \pi(u) p(u,v) (f(u) - f(v))^2 \\ &= \inf_{f \in \mathcal{F}} \sum_{[u,v] \in E(\mathbb{H}')} \frac{2}{D+1} (f(u) - f(v))^2 \\ &= \inf_{f \in \mathcal{F}} \sum_{[u,v] \in E(\mathbb{H})} \frac{2}{D+1} (f(u) - f(v))^2. \end{aligned}$$

In particular, let $\{Y_n\}$ be the simple random walk in $\mathbb{G}$ and

$$\sigma_N = \min\{n \geq 0; |Y_n|_\mathbb{G} = N\},$$
$$\sigma_o^+ = \min\{n \geq 1; Y_n = o\},$$
$$\rho(m) = P(\sigma_N < \sigma_o^+ | Y_0 = m).$$

Then $f(x) = \rho(m(x))$ is in $\mathcal{F}$. Plugging it into (3.5), in light of (3.3), we conclude that

$$\sum_{\substack{x \in V(\mathbb{H}) \\ m(x)=o}} P_\omega(\tau_N < \tau_o^+ | Z(0) = x) \leq \frac{1}{D+1} \sum_{[u,v] \in E(\mathbb{H})} (\rho(m(u)) - \rho(m(v)))^2.$$

Note that $P_\omega(\tau_N < \tau_o^+ | Z(0) = x) = 0$ if $m(x) = o$ and there is no $y$ in the cluster such that $|m(y)|_\mathbb{G} = N$. Summing over all disjoint clusters, we get

$$\sum_{x : m(x)=o} P_\omega(\tau_N < \tau_o^+ | Z(0) = x)$$

$$= \sum_\mathbb{H} \sum_{x \in V(\mathbb{H}), m(x)=o} P_\omega(\tau_N < \tau_o^+ | Z(0) = x)$$

$$\leq \frac{|\mathbb{F}|^{|\mathbb{G}_N|}}{D+1} \sum_{[u,v] \in E(\mathbb{G}_N)} (\rho(u) - \rho(v))^2.$$

Averaging over realizations of percolation in $\mathbb{G}_N \ltimes \sum_{x \in \mathbb{G}_N} \mathbb{F}$, we see that $\mathbf{E}P_\omega(\tau_N < \tau_o^+ | Z(0) = (o, \eta))$ is independent of $\eta$. There are $|\mathbb{F}|^{|\mathbb{G}_N|}$ vertices in $\mathbb{G}_N \ltimes \sum_{x \in V(\mathbb{G}_N)} \mathbb{F}$ with the marker at $o$. Therefore,

$$|\mathbb{F}|^{|\mathbb{G}_N|} \mathbf{E}P_\omega(\tau_N < \tau_o^+ | Z(0) = (o, \eta))$$

$$\leq \frac{|\mathbb{F}|^{|\mathbb{G}_N|}}{D+1} \sum_{[u,v] \in E(\mathbb{G}_N)} (\rho(u) - \rho(v))^2.$$

After cancellation,

$$\mathbf{E}P_\omega(\tau_N < \tau_o^+ | Z(0) = (o, \eta))$$

$$\leq \frac{1}{D+1} \sum_{[u,v] \in E(\mathbb{G}_N)} (\rho(u) - \rho(v))^2$$

$$= \frac{1}{D+1} P(\sigma_N < \sigma_o^+ | Y_0 = o) \to 0 \quad \text{as } N \to \infty,$$

since the simple random walk in $\mathbb{G}$ is recurrent. □

PROOF OF THEOREM 3.1. Let $Z$ be the delayed simple random walk in the infinite cluster. By Lemma 3.3, it suffices to show that $\lim_n \|Z(n)\|/n = 0$ a.s.

Let $\mathcal{R}_n = \{m_0, m_1, \ldots, m_n\} \subset V(\mathbb{G})$ be the range of the marker up to time $n$. Then

$$|\mathcal{R}_n| = 1 + \sum_{k=0}^{n-1} \mathbb{1}_{\{m_k \neq m_{k+1}, m_k \neq m_{k+2}, \ldots, m_k \neq m_n\}}$$

(3.6)
$$\leq \ell + \sum_{k=0}^{n-\ell} \mathbb{1}_{\{m_k \neq m_{k+1}, m_k \neq m_{k+2}, \ldots, m_k \neq m_{k+\ell}\}}$$

for any fixed integer $\ell$. As explained in Lyons and Schramm (1999), in the large probability space where both the percolation and the walk $Z$ are defined, the law of the infinite cluster $\omega$ as seen from the walker $Z(n)$ is stationary. Therefore

$$\{\mathbb{1}_{\{m_k \neq m_{k+1}, m_k \neq m_{k+2}, \ldots, m_k \neq m_{k+\ell}\}}; k = 0, 1, 2, 3, \ldots\}$$

is a stationary sequence in the large space. By (3.6) and the Birkhoff ergodic theorem,

$$\mathbf{E}E_\omega \limsup_n \frac{|\mathcal{R}_n|}{n} \leq \mathbf{E}E_\omega \lim_n \frac{1}{n} \sum_{k=0}^{n-\ell} \mathbb{1}_{\{m_k \neq m_{k+1}, m_k \neq m_{k+2}, \ldots, m_k \neq m_{k+\ell}\}}$$
$$= \mathbf{E}P_\omega(m_0 \neq m_1, m_0 \neq m_2, \ldots, m_0 \neq m_\ell).$$

By Lemma 3.4, the right-hand side tends to 0 as $\ell \to \infty$. Consequently, for a.e. $\omega$ and $Z$,

$$\lim_n \frac{|\mathcal{R}_n|}{n} = 0.$$

Note that $\mathcal{R}_n$ is connected in $V(\mathbb{G})$, and all sites in $\mathcal{R}_n$ can be visited within at most $2|\mathcal{R}_n|$ steps using depth-first search along a spanning tree in $\mathcal{R}_n$. Thus in $\mathbb{W}$,

$$\|Z(n)\| \leq |m_n|_\mathbb{G} + 2|\mathcal{R}_n| + \sum_{x \in \mathcal{R}_n} |\eta_n(x)|_\mathbb{F} \leq (1 + 2 + |\mathbb{F}|)|\mathcal{R}_n|.$$

We conclude that $\mathbf{E}E_\omega \limsup_{n \to \infty} \|Z(n)\|/n = 0$. □

Our next goal is to prove Theorem 3.2.

LEMMA 3.5. *Suppose that the infinite cluster of $p$-Bernoulli bond percolation on the Cayley graph $\mathbb{G}$ is transient. Let $Z$ be the delayed random walk in the infinite cluster of $p$-Bernoulli bond percolation on $\mathbb{W} = \mathbb{G} \ltimes \sum_{x \in \mathbb{G}} \mathbb{F}$ and write $Z(n) = (m_n, \eta_n)$. Then*

$$\mathbf{E}P_\omega(m_n \neq m_0 \text{ for all } n \geq 1) > 0.$$

PROOF. We shall prove that

(3.7) $$\lim_{N \to \infty} \mathbf{E}P_\omega(\tau_N < \tau_o^+ | m_0 = o) > 0,$$

where $\tau_N$ and $\tau_o^+$ are stopping times defined in the proof of Lemma 3.4.

Recall the finite graphs $\mathbb{G}_N$ and $\mathbb{G}_N \ltimes \sum_{x \in \mathbb{G}_N} \mathbb{F}$ defined in the proof of Lemma 3.4. Vertices $(m, \eta)$ of $\mathbb{G}_N \ltimes \sum_{x \in \mathbb{G}_N} \mathbb{F}$ are classified into $|\mathbb{F}|^{|\mathbb{G}_N|}$ classes according to the second component $\eta$. For a fixed configuration $\eta$, denote by $\mathbb{G}_N(\eta)$ the subgraph induced by the class of vertices $\{(m, \eta); m \in V(\mathbb{G}_N)\}$. Clearly, there is a graph isomorphism between $\mathbb{G}_N(\eta)$ and $\mathbb{G}_N$ for any $\eta \in \sum_{x \in \mathbb{G}_N} \mathbb{F}$. Let the cluster within $\mathbb{G}_N(\eta)$ containing $(o, \eta)$ be

$$C_o(\eta) = \{(m, \eta); (m, \eta) \leftrightarrow (o, \eta) \text{ within } \mathbb{G}_N(\eta)\}.$$

Run a simple random walk $\{(Y_j^\eta, \eta)\}_{j \geq 0}$ in $C_o(\eta)$ starting from $(o, \eta)$. Recall $\sigma_N = \min\{j \geq 0; |Y_j^\eta|_{\mathbb{G}} = N\}$ and $\sigma_o^+ = \min\{j \geq 1; Y_j^\eta = o\}$. Then $P(\sigma_N < \sigma_o^+ | Y_0^\eta = o)$ is decreasing in $N$. The hypothesis of the lemma (transience of the infinite cluster) means that

(3.8) $$\lim_{N \to \infty} \mathbf{E} P_\omega(\sigma_N < \sigma_o^+ | Y_0^\eta = o) > 0.$$

There may be several disjoint clusters in a realization of $p$-Bernoulli bond percolation in $\mathbb{G}_N \ltimes \sum_{x \in \mathbb{G}_N} \mathbb{F}$ and each cluster may have several vertices with the marker at $o$. Take a cluster, say $\mathbb{H}$, and run the delayed simple random walk $Z$ in $V(\mathbb{H})$. It follows from (3.3) and (3.5) that

(3.9)
$$\sum_{x \in V(\mathbb{H}), m(x) = o} P_\omega(\tau_N < \tau_o^+ | Z(0) = x)$$
$$= \inf_{f \in \mathcal{F}} \sum_{[u,v] \in E(\mathbb{H})} \frac{1}{d_\mathbb{G} + d_\mathbb{F} + 1} (f(u) - f(v))^2,$$

where $\mathcal{F}$ is the class of functions satisfying (3.4), and $d_\mathbb{G}$ and $d_\mathbb{F}$ are the degrees of a vertex of $\mathbb{G}$ and $\mathbb{F}$, respectively. Notice that (3.9) is still valid even if there is no vertex $y \in \mathbb{H}$ such that $|m(y)|_G = N$. Summing over all disjoint clusters, we get

(3.10)
$$\sum_\eta P_\omega(\tau_N < \tau_o^+ | Z(0) = (o, \eta))$$
$$= \sum_\mathbb{H} \sum_{(o,\eta) \in V(\mathbb{H})} P_\omega(\tau_N < \tau_o^+ | Z(0) = (o, \eta))$$
$$= \inf_{f \in \mathcal{F}} \sum \frac{1}{d_\mathbb{G} + d_\mathbb{F} + 1} (f(u) - f(v))^2,$$

where the summation is over all open bonds $[u, v]$ of $\mathbb{G}_N \ltimes \sum_{x \in \mathbb{G}_N} \mathbb{F}$. Every term in (3.10) is nonnegative. Discarding those terms involving an open edge $[u, v]$ where $u$ and $v$ are in different classes (i.e., the markers of $u$ and $v$ are the same), we get the following inequality:

$$\text{RHS of (3.10)} \geq \inf_{f \in \mathcal{F}} \sum_\eta \sum_{[u,v] \in E(\mathbb{G}_N(\eta)), \text{ open}} \frac{1}{d_\mathbb{G} + d_\mathbb{F} + 1} (f(u) - f(v))^2$$

(3.11)
$$= \frac{d_{\mathbb{G}}+1}{d_{\mathbb{G}}+d_{\mathbb{F}}+1} \sum_{\eta} P_\omega(\sigma_N < \sigma_o^+ | Y_0^\eta = o).$$

Combining (3.11) with (3.10), we conclude that for a.e. realization $\omega$ of the Bernoulli bond percolation,

$$(3.12) \quad \sum_{\eta} P_\omega(\tau_N < \tau_o^+ | Z(0) = (o,\eta)) \geq \sum_{\eta} \frac{d_{\mathbb{G}}+1}{d_{\mathbb{G}}+d_{\mathbb{F}}+1} P_\omega(\sigma_N < \sigma_o^+ | Y_0^\eta = o).$$

Taking expectation over the Bernoulli bond percolation $\omega$, we find that $\mathbf{E}P_\omega(\tau_N < \tau_o^+ | Z(0) = (o,\eta))$ is independent of $\eta$. It follows from (3.12) that

$$\mathbf{E}P_\omega(\tau_N < \tau_o^+ | Z(0) = (o,\eta)) \geq \frac{d_{\mathbb{G}}+1}{d_{\mathbb{G}}+d_{\mathbb{F}}+1} \mathbf{E}P_\omega(\sigma_N < \sigma_o^+ | Y_0^\eta = o).$$

Taking the limit as $N \to \infty$, inequality (3.7) then follows from (3.8). $\square$

PROOF OF THEOREM 3.2. The delayed simple random walk $Z = Z^\omega$ on the infinite cluster $\omega$ is a reversible Markov chain with respect to the uniform measure on $\omega$. In conjunction with the stationarity of $Z$ in the big space, this gives

$$(3.13) \quad \mathbf{E}P_\omega(m_i \neq m_n, 0 \leq i \leq n-1) = \mathbf{E}P_\omega(m_i \neq m_0, 1 \leq i \leq n).$$

Define
$$\zeta(k) = \begin{cases} 1, & \text{if } \eta_k \neq \eta_{k-1}, m_i \neq m_k \text{ for } 0 \leq i \leq k-2 \text{ and for } i \geq k+1, \\ 0, & \text{otherwise.} \end{cases}$$

Then
$$E_\omega \zeta(k) = P_\omega(\eta_k \neq \eta_{k-1} \text{ and } m_i \neq m_k \text{ for } 0 \leq i \leq k-2 \text{ and for } i \geq k+1)$$

can be written as a product of three terms:

$$P_\omega(m_i \neq m_k \text{ for all } i \geq k+1 | \eta_k \neq \eta_{k-1} \text{ and } m_j \neq m_k \text{ for } 0 \leq j \leq k-2)$$
$$\times P_\omega(\eta_k \neq \eta_{k-1} | m_j \neq m_{k-1} \text{ for } 0 \leq j \leq k-2)$$
$$\times P_\omega(m_j \neq m_{k-1} \text{ for } j = 0, 1, 2, \ldots, k-2).$$

In the big probability space, the distribution of $(Z, \omega)$ is invariant under the shift by Lemma 4.1 of BLS (1999). Taking the expectation over realizations of $\omega$, then

$$\mathbf{E}P_\omega(m_i \neq m_k \text{ for all } i \geq k+1 | \eta_k \neq \eta_{k-1} \text{ and } m_j \neq m_k \text{ for } 0 \leq j \leq k-2)$$
$$= \mathbf{E}P_\omega(m_i \neq m_0 \text{ for all } i \geq 1).$$

The last equality holds because when a light at location $m_k$ is switched for the first time at step $k$, and then the marker moves away from that site, a regeneration occurs: while $m_i \neq m_k$, the walker $Z(\cdot)$ is traveling in virgin territory that was not explored prior to time $k$.

Moreover,

$$\mathbf{E} P_\omega(\eta_k \neq \eta_{k-1} | m_j \neq m_{k-1} \text{ for } 0 \leq j \leq k-2) \geq \frac{pd_\mathbb{F}}{d_\mathbb{G} + d_\mathbb{F} + 1};$$

and using reversibility (3.13),

$$\mathbf{E} P_\omega(m_i \neq m_{k-1} \text{ for } 0 \leq i \leq k-2) \geq \mathbf{E} P_\omega(m_i \neq m_0 \text{ for all } i \geq 1).$$

Therefore

$$\mathbf{E} E_\omega \zeta(k) \geq \frac{pd_\mathbb{F}}{d_\mathbb{G} + d_\mathbb{F} + 1} (\mathbf{E} P_\omega(m_i \neq m_0 \text{ for all } i \geq 1))^2.$$

Finally, because $\|Z(n)\| \geq \sum_{k=1}^n \zeta(k)$,

$$\mathbf{E} E_\omega \lim_n \frac{\|Z(n)\|}{n} = \lim_n \mathbf{E} E_\omega \frac{\|Z(n)\|}{n}$$
$$\geq \lim_n \mathbf{E} E_\omega \frac{1}{n} \sum_{k=1}^n \zeta(k) \geq \lim_n \frac{1}{n} \sum_{k=1}^n \mathbf{E} E_\omega \zeta(k) > 0.$$

Since $\lim_n \|Z(n)\|/n$ exists and is a constant a.s., it must be positive, and we are done. $\square$

## APPENDIX

The goal of this Appendix is to prove the following sharpening of Theorem 1.1.

THEOREM A.1. *Consider $p$-Bernoulli bond percolation on a graph $G$ with $\iota_E^*(G) > 0$. If $p > 1/(1 + \iota_E^*(\mathbb{G}))$, then almost surely on the event that the open cluster $\mathbb{H}$ containing $o$ is infinite, it satisfies $\iota_E^*(\mathbb{H}) > 0$.*

PROOF. We will use the notation and some of the ideas of the proof of Theorem 1.1. First note that in $p$-Bernoulli bond percolation, for any $0 < \alpha < p$, we can estimate the conditional probability

$$(A.1) \quad P\left(\frac{|\partial_\mathbb{H} S|}{|\partial_\mathbb{G} S|} \leq \alpha \Big| S \in \mathcal{A}_n(\mathbb{H})\right) = P(\text{Binom}(n,p) \leq \alpha n) \leq e^{-n I_p(\alpha)},$$

where the rate function $I_p(\cdot)$ is continuous, and

$$-\log(1-p) = I_p(0) > I_p(\alpha) > 0 \qquad \text{for } 0 < \alpha < p.$$

Therefore,

$$P\left(\exists S \in \mathcal{A}_n(\mathbb{H}): \frac{|\partial_\mathbb{H} S|}{|\partial_\mathbb{G} S|} \leq \alpha\right)$$

$$\leq \sum_{S \in \mathcal{A}_n} P\left(S \in \mathcal{A}_n(\mathbb{H}), \frac{|\partial_\mathbb{H} S|}{|\partial_\mathbb{G} S|} \leq \alpha\right)$$

(A.2)
$$\leq \sum_{S \in \mathcal{A}_n} e^{-nI_p(\alpha)} P(S \in \mathcal{A}_n(\mathbb{H}))$$

$$= e^{n(I_p(0)-I_p(\alpha))} \sum_{S \in \mathcal{A}_n} (1-p)^n P(S \in \mathcal{A}_n(\mathbb{H}))$$

$$= e^{n(I_p(0)-I_p(\alpha))} P(|V(\mathbb{H})| < \infty, |\partial_\mathbb{G} V(\mathbb{H})| = n),$$

where the last step used the identity $(1-p)^n P(S \in \mathcal{A}_n(\mathbb{H})) = P(\mathbb{H} = S)$ for $S \in \mathcal{A}_n$.

To estimate $P(|V(\mathbb{H})| < \infty, |\partial_\mathbb{G} V(\mathbb{H})| = n)$, we use the method of Theorem 2 of Benjamini and Schramm (1996); see also Theorem 6.18 in Lyons and Peres (2004). Let us briefly recall that argument. Choose $h < \iota_E^*(\mathbb{G})$ such that $p > \frac{1}{1+h}$. Then there exists $n_h < \infty$ such that $|\partial_\mathbb{G} S|/|S| > h$ for all $S \in \mathcal{A}_n$ with $n > n_h$. Fix an ordering of the edges $E(\mathbb{G}) = (e_1, e_2, \ldots)$ such that $o$ is an endpoint of $e_1$, and take two i.i.d. sequences $\{Y_i\}$ and $\{Y_i'\}$ of Bernoulli$(p)$ variables. Build recursively the percolation cluster $\mathbb{H}$ of $o$, together with its boundary $\partial_\mathbb{G} V(\mathbb{H})$, using the sequence $\{Y_i\}$, as follows. At step zero, we start with $\mathbb{H}_0$ consisting just of $o$. In step $j \geq 1$, consider the first unexamined edge $e_{n_j}$ in the ordering above that has one endpoint in $V(\mathbb{H}_{j-1})$, and one endpoint in its complement. (If there is no such edge, the process stops and we have $\mathbb{H} = \mathbb{H}_{j-1}$.) Let $\mathbb{H}_j$ be $\mathbb{H}_{j-1}$ with $e_{n_j}$ added if $Y_j = 1$, and $\mathbb{H}_j = \mathbb{H}_{j-1}$ if $Y_j = 0$. If the process continues indefinitely, then the increasing union of all the $\mathbb{H}_j$ is the infinite cluster $\mathbb{H}$. Having finished with growing the finite or infinite cluster $\mathbb{H}$, build the remainder of the percolation configuration using the sequence $\{Y_i'\}$.

If the process terminates after $N$ steps with a finite cluster $\mathbb{H} = \mathbb{H}_N$ that has $v$ vertices and $n$ closed boundary edges, then $N \geq n + v - 1$ and $\sum_{j=1}^N Y_j = v - 1$. For $n > n_h$ we must have $n > vh$, whence $v - 1 < N/(1+h)$. It follows that

$$\{|V(\mathbb{H})| < \infty, \, |\partial_\mathbb{G} V(\mathbb{H})| = n\} \subset \bigcup_{N=n}^\infty B_N,$$

where

$$B_N = \left\{\sum_{j=1}^N Y_j \leq \frac{N}{1+h}\right\}.$$

By the large deviation principle, $P(B_N) \leq e^{-N\delta_p}$, where $\delta_p = I_p(\frac{1}{1+h}) > 0$, since $p > \frac{1}{1+h}$. Thus for some constant $C_p < \infty$,

$$(A.3) \qquad P(|V(\mathbb{H})| < \infty, |\partial_\mathbb{G} V(\mathbb{H})| = n) \leq \sum_{N=n}^{\infty} e^{-N\delta_p} \leq C_p e^{-n\delta_p}.$$

Taking $\alpha > 0$ in (A.2) so small that $I_p(0) - I_p(\alpha) < \delta_p$, we deduce that (A.2) is summable in $n$. An application of the Borel–Cantelli lemma, just as at the end of the proof of Theorem 1.1, gives that

$$\imath_E^*(\mathbb{H}) \geq \alpha \imath_E^*(\mathbb{G}) > 0$$

almost surely on the event that $\mathbb{H}$ is infinite. $\square$

For site percolation on $\mathbb{G}$, the vertex version of anchored expansion is the relevant notion. Let $\partial^V S$ denote the set of vertices in $S^c$ having a neighbor in $S$, and suppose that

$$\imath_V^*(\mathbb{G}) := \lim_{n \to \infty} \inf \left\{ \frac{|\partial^V S|}{|S|} : o \in S \subset V(\mathbb{G}), S \text{ is connected, } n \leq |S| < \infty \right\} > 0.$$

Then the corresponding form of (A.2) needs no modification, while the analog of (A.3) can be proved using an ordering of the vertices $V(\mathbb{G})$. Hence, the following result holds:

THEOREM A.2. *Consider $p$-Bernoulli site percolation on a graph $G$ with $\imath_V^*(G) > 0$. If $p > 1/(\imath_V^*(\mathbb{G})+1)$, then almost surely on the event that the open cluster $\mathbb{H}$ containing $o$ is infinite, it satisfies $\imath_V^*(\mathbb{H}) > 0$.*

These results are sharp for the $(b+1)$-regular trees $\mathbb{T}_b$, for which $p_c(\mathbb{T}_b) = 1/b$ for both bond and site percolations, while $\imath_E^*(\mathbb{T}_b) = \imath_V^*(\mathbb{T}_b) = b - 1$.

**Acknowledgments.** We are grateful to Noam Berger, Alan Hammond, David Revelle, Bálint Virág and the referee for helpful comments.

## REFERENCES

ANGEL, O., BENJAMINI, I., BERGER, N. and PERES, Y. (2004). Transience of percolation clusters on wedges. *Electron. J. Probab.* To appear.

BENJAMINI, I., LYONS, R. and SCHRAMM, O. (1999). Percolation perturbations in potential theory and random walks. In *Random Walks and Discrete Potential Theory* (M. Picardello and W. Woess, eds.) 56–84. Cambridge Univ. Press. MR1802426

BENJAMINI, I. and SCHRAMM, O. (1996). Percolation beyond $\mathbb{Z}^d$, many questions and a few answers. *Electron. Comm. Probab.* **1** 71–82. MR1423907

BENJAMINI, I. and SCHRAMM, O. (1998). Oriented random walk on the Heisenberg group and percolation. Unpublished manuscript.

ANCHORED EXPANSION, PERCOLATION AND SPEED 19Burton, R. M. and Keane, M. (1989). Density and uniqueness in percolation. *Comm. Math. Phys.* **121** 501–505. [MR990777](MR990777)

Dembo, A. and Zeitouni, O. (1998). *Large Deviations Techniques and Applications*, 2nd ed. Springer, New York. [MR1619036](MR1619036)

Grimmett, G. R., Kesten, H. and Zhang, Y. (1993). Random walk on the infinite cluster of the percolation model. *Probab. Theory Related Fields* **96** 33–44. [MR1222363](MR1222363)

Häggström, O., Schonmann, R. and Steif, J. (2000). The Ising model on diluted graph and strong amenability. *Ann. Probab.* **28** 1111–1137. [MR1797305](MR1797305)

Harris, T. E. (1963). *The Theory of Branching Processes*. Springer, Berlin. [MR163361](MR163361)

Kaimanovich, V. A. and Vershik, A. M. (1983). Random walks on discrete groups: Boundary and entropy. *Ann. Probab.* **11** 457–490. [MR704539](MR704539)

Kesten, H. (1982). *Percolation Theory for Mathematicians*. Birkhäuser, Boston. [MR692943](MR692943)

Liggett, T. M. (1985). *Interacting Particle Systems.* Springer, New York. [MR776231](MR776231)

Lyons, R. (1992). Random walks, capacity, and percolation on trees. *Ann. Probab.* **20** 2043–2088. [MR1188053](MR1188053)

Lyons, R., Pemantle, R. and Peres, Y. (1995). Ergodic theory on Galton–Watson trees: Speed of random walk and dimension of harmonic measure. *Ergodic Theory Dynam. Systems* **15** 593–619. [MR1336708](MR1336708)

Lyons, R., Pemantle, R. and Peres, Y. (1996). Random walks on the lamplighter group. *Ann. Probab.* **24** 1993–2006. [MR1415237](MR1415237)

Lyons, R. and Peres, Y. (2004). *Probability on Trees and Networks*. Cambridge Univ. Press. To appear. Available at http://mypage.iu.edu/˜rdlyons.

Lyons, R. and Schramm, O. (1999). Indistinguishability of percolation clusters. *Ann. Probab.* **27** 1809–1836. [MR1742889](MR1742889)

Thomassen, C. (1992). Isoperimetric inequalities and transient random walks on graphs. *Ann. Probab.* **20** 1592–1600. [MR1175279](MR1175279)

Varopoulos, N. Th. (1985). Long range estimates for Markov chains. *Bull. Sci. Math. (2)* **109** 225–252. [MR822826](MR822826)

Virág, B. (2000). Anchored expansion and random walk. *Geom. Funct. Anal.* **10** 1588–1605. [MR1810755](MR1810755)

Woess, W. (2000). *Random Walks on Infinite Graphs and Groups*. Cambridge Univ. Press. [MR1743100](MR1743100)
School of Mathematical Sciences
Peking University
Beijing 100871
China
e-mail: dayue@math.pku.edu.cn

Department of Statistics
University of California
Berkeley, California 94720
USA
e-mail: peres@stat.berkeley.edu
e-mail: gabor@stat.berkeley.edu